\documentclass[12pt]{article}

\usepackage{amsmath,amsfonts,amssymb,amsthm,cite}
\usepackage{geometry}

\theoremstyle{plain}
\newtheorem{theorem}{Theorem}

\newtheorem{corollary}{Corollary}
\newtheorem{definition}{Definition}
\theoremstyle{definition}

\newcommand{\nonprint}[1]{}

\begin{document}

\begin{flushleft}

\vspace{+0.3cm}
\large

\begin{center}
\textbf{On Differential Systems in Sobolev spaces with Generic Inhomogeneous Boundary Conditions}\\
\end{center}

\begin{center}
\vspace{+0.4cm}
\textbf{Olena Atlasiuk, Vladimir Mikhailets}
\vspace{+0.3cm}
\end{center}

\end{flushleft}

\normalsize

\begin{abstract}
We study linear systems of ordinary differential equations of an arbitrary order on a finite interval with the most general (generic) inhomogeneous boundary conditions in Sobolev spaces. We investigate the character of solvability of inhomogeneous boundary-value problems, prove their Fredholm properties, and find the indices, the dimensions of the kernel, and the cokernel of these problems. Moreover, we obtained necessary and sufficient conditions for continuity in a parameter of solutions to the introduced problems in Sobolev spaces.
\end{abstract}



\section{Boundary-value problems}
Boundary-value problems with inhomogeneous boundary conditions containing derivatives whose order is greater than or equal to the order of the differential equation arise naturally in some mathematical models (see, for example, \cite{Vent1959,Luo1991,Kra1961}). The theory of such problems contains few results so far even for the case of ordinary differential equations. The purpose of this paper is to develop this theory regarding linear systems of ordinary differential equations of an arbitrary order with generic inhomogeneous boundary conditions in Sobolev spaces.

Let a finite interval  $(a,b)\subset\mathbb{R}$ and the next parameters be given $$\{m, \, n+1, \,r,\,l\} \subset \mathbb{N}, \, 1\leq p\leq \infty.$$ By $W_p^{n+r}=W_p^{n+r}\bigl([a,b];\mathbb{C}\bigr):= \bigl\{y\in C^{n+r-1}[a,b] \colon y^{(n+r-1)}\in AC[a,b], \, y^{(n+r)}\in L_p[a,b]\bigr\}$ we denote a complex Sobolev space and set $W_p^0:=L_p$. This space is a Banach one with respect to the norm
$$
\bigl\|y\bigr\|_{n+r,p}=\sum_{k=0}^{n+r}\bigl\|y^{(k)}\bigr\|_{p},
$$
where $\|\cdot\|_p$ is the norm in $L_p\bigl([a,b]; \mathbb{C}\bigr)$. Similarly, by $(W_p^{n+r})^{m}:=W_p^{n+r}\bigl([a,b];\mathbb{C}^{m}\bigr)$ and $(W_p^{n+r})^{m\times m}:=W_p^{n+r}\bigl([a,b];\mathbb{C}^{m\times m}\bigr)$
we denote Sobolev spaces of vector-valued functions and matrix-valued functions, respectively, whose ele\-ments belong to the function space $W_p^{n+r}$.

We consider the following linear boundary-value problem
\begin{equation}\label{eq1}
    (Ly)(t):=y^{(r)}(t) + \sum\limits_{j=1}^rA_{r-j}(t)y^{(r-j)}(t)=f(t),\quad
t \in(a,b),
\end{equation}
\begin{equation}\label{eq2}
    By= c,
\end{equation}
\noindent where matrix-valued functions $A_{r-j}(\cdot)\in (W_{p}^{n})^{m\times m}$, a vector-valued function $f(\cdot)\in (W_{p}^{n})^{m}$, vector $c\in\mathbb{C}^{l}$, a linear continuous operator
\begin{equation}\label{eq3}
B\colon(W_{p}^{n+r})^{m} \rightarrow\mathbb{C}^{l}
\end{equation}
are arbitrarily chosen; and a vector-valued function $y(\cdot)\in (W_{p}^{n+r})^m$ is unknown. If $l<r$, then the boundary conditions are underdetermined. If $l>r$, then the boundary conditions are overdetermined.

We represent vectors and vector-valued functions in the form of columns. A solution to the problem \eqref{eq1}, \eqref{eq2} is understood as a vector-valued function $y(\cdot)\in (W_{p}^{n+r})^m$ satisfying equation \eqref{eq1} (for $n\geq 1$ everywhere, and for $n= 0$ almost everywhere) on $(a,b)$, and equality \eqref{eq2} specifying $l$ scalar boundary conditions.

It includes all known types of classical boundary conditions, namely, the Cauchy problem, two- and many-point problems, integral and mixed problems, and numerous nonclassical problems. The last class of problems may contain derivatives (generally fractional) $y^{(k)}(\cdot)$, with $0 < k \leq n+r$, (see, for instance, \cite{Kilbas2006}).

The solutions of equation \eqref{eq1} fill the space $(W_{p}^{n+r})^m$ if its right-hand side $f(\cdot)$ runs through the space $(W_{p}^{n})^m$. Therefore, the boundary condition \eqref{eq2} with continuous operator \eqref{eq3} is the most general condition for this equation.

For $1\leq p < \infty$, every operator $B$ in \eqref{eq3}
admits a unique analytic representation
$$
By=\sum _{k=0}^{n+r-1} \alpha_{k} y^{(k)}(a)+\int_{a}^b \Phi(t)y^{(n+r)}(t){\rm d}t, \quad y(\cdot)\in (W_{p}^{n+r})^{m},
$$
where the matrices $\alpha_{k}\in\mathbb{C}^{rm\times m}$, \smash{$1/p + 1/p^{'}=1$}, and the matrix-valued function $\Phi(\cdot)\in L_{p^{'}}\bigl([a, b]; \mathbb{C}^{rm\times m}\bigr)$.

For $p=\infty$ this formula also defines a bounded operator $B\colon (W_{\infty}^{n+r})^{m} \rightarrow \mathbb{C}^{rm}$. However, there exist other operators from this class generated by the integrals over finitely additive measures. Hence, the study of the $p=\infty$ case faces additional difficulties, unlike when $p< \infty$ (\cite{KM2013,GKM2017,Atl2}).

\section{Solvability}

With the generic inhomogeneous boundary-value problem \eqref{eq1}, \eqref{eq2}, we associate a linear operator in pair of Banach spaces
\begin{equation}\label{eq4}
(L,B)\colon (W^{n+r}_p)^m\rightarrow (W^{n}_p)^m\times\mathbb{C}^l.
\end{equation}

Recall that a linear continuous operator $T\colon X \rightarrow Y$, where $X$ and $Y$ are Banach spaces, is called a Fredholm
operator if its kernel and cokernel are finite-dimensional. If operator $T$ is Fredholm, then its
range $T(X)$ is closed in $Y$ and the index is finite
$$
\mathrm{ind}\,T:=\dim\ker T-\dim(Y/T(X))\in \mathbb{Z}.
$$

\begin{theorem}\label{th1} The linear operator \eqref{eq4} is a bounded Fredholm operator with index $mr-l$.
\end{theorem}

The proof of Theorem \ref{th1} uses the well-known theorem on the stability of the index of a linear operator with respect to compact additive perturbations (see \cite{Atl1}).

This theorem naturally raises the question of finding the Fredholm numbers (i.e. the dimensions of the problem kernel and co-kernel). This is a quite difficult task because the Fredholm numbers may vary even under arbitrarily small one-dimensional perturbations.

To formulate the following result, let us introduce some notation and definitions.

For each number $k \in \{1,\dots, r\}$, we consider a family of matrix Cauchy problems:
$$
Y_k^{(r)}(t)+\sum\limits_{j=1}^rA_{r-j}(t)Y_k^{(r-j)}(t)=O_{m},\quad t\in (a,b),
$$
with the initial conditions
$$
Y_k^{(j-1)}(a) = \delta_{k,j}I_m,\qquad j \in \{1,\dots, r\}.
$$

Here, $Y_k(\cdot)$ is an unknown $(m \times m)$ -- matrix-valued function, and $\delta_{k,j}$ is the Kronecker symbol.

By $\left[BY_k\right]$ we denote the numerical $(m\times l)$ -- matrix, in which $j$-th column is the result of action of the operator $B$ on the $j$-th column of the matrix-valued function $Y_k(\cdot)$.

\begin{definition}\label{d1}  A block rectangular numerical matrix
$$
M(L,B):=\left(\left[BY_0\right],\dots,\left[BY_{r-1}\right]\right) \in \mathbb{C}^{l\times mr}
$$
is characteristic to the inhomogeneous boundary-value problem \eqref{eq1}, \eqref{eq2}.
\end{definition}

It consists of $r$ rectangular block columns $\left[BY_k\right]\in \mathbb{C}^{m\times l}$.

Here $mr$ is the number of scalar differential equations of the system \eqref{eq1}, and $l$ is the number of scalar boundary conditions.

\begin{theorem}\label{th2} The dimensions of the kernel and cokernel of the operator \eqref{eq4} are equal to the~dimensions of the kernel and cokernel of the characteristic matrix $M(L,B)$, respectively.
\end{theorem}

Theorem \ref{th2} implies necessary and sufficient conditions for the invertibility of the operator \eqref{eq4}.

\begin{corollary}\label{co1} The operator $(L,B)$ is invertible if and only if $l=mr$ and the square matrix $M(L,B)$ is nondegenerate.
\end{corollary}

The following theorem shows that the functions $\operatorname{dim} \operatorname{ker}(L,B)$ and $\operatorname{dim} \operatorname{coker}(L,B)$ are semi-continuous in the strong operator topology in the class of problems we have considered.

With the problem \eqref{eq1}, \eqref{eq2}, we consider a sequence of boundary-value problems
\begin{equation}\label{eq5}
L(k)y(t,k):= y^{(r)}(t,k)+\sum_{j=1}^{r}A_{r-j}(t,k)y^{(r-j)}(t,k)=f(t,k),\quad t\in (a,b),
\end{equation}
\begin{equation}\label{eq6}
B(k)y(\cdot,k)=c(k), \quad k\in\mathbb{N},
\end{equation}
where the matrix-valued functions $A_{r-j}(\cdot,k)$, the vector-valued function $f(\cdot,k)$, the vector $c(k)$, and the linear continuous operators $B(k)$ satisfy the above conditions to the problem \eqref{eq1}, \eqref{eq2}.

With the boundary-value problem \eqref{eq5}, \eqref{eq6}, we associate a sequence of linear continuous operators
$$
(L(k),B(k))\colon(W^{n+r}_p)^{m}\rightarrow (W^{n}_p)^{m}\times\mathbb{C}^{l}
$$
and a sequence of characteristic matrices depending on the parameter $k\in\mathbb{N}$
$$
M\big(L(k),B(k)\big):=\big(\left[B(k)Y_0(k)\right],\dots,\left[B(k)Y_{r-1}(k)\right]\big) \subset \mathbb{C}^{mr\times l}.
$$

We now formulate a sufficient condition for the convergence of the characteristic matrices $M\left(L(k),B(k)\right)$ to the matrix $M\left(L,B\right)$.

\begin{theorem}\label{th3} If the sequence of operators $\big(L(k),B(k)\big)$ converges strongly to the operator $\big(L,B\big)$ for $k\rightarrow \infty$, then the sequence of characteristic matrices $M\big(L(k),B(k)\big)$ converges to the matrix $M\big(L,B\big)$.
\end{theorem}


\begin{corollary}\label{co2} Under the assumptions from Theorem \ref{th3}, the following inequalities hold for sufficiently large $k$:
\begin{gather*}
\operatorname{dim} \operatorname{ker}\left(L(k),B(k)\right)\leq \operatorname{dim} \operatorname{ker}\left(L,B\right),  \\
\operatorname{dim} \operatorname{coker}\left(L(k),B(k)\right)\leq \operatorname{dim} \operatorname{coker}\left(L,B\right).
\end{gather*}
\end{corollary}

In particular:
\begin{itemize}
\item[1.] If $l=mr$ and the operator $(L,B)$ is invertible, then the operators $\left(L(k),B(k)\right)$ are also invertible for large $k$.

\item[2.] If the boundary-value problem \eqref{eq1}, \eqref{eq2} has a solution for any values of the right-hand sides, then the boundary-value problems \eqref{eq5}, \eqref{eq6} also have the solutions for large $k$.

\item[3.] If the boundary-value problem \eqref{eq1}, \eqref{eq2} has a unique solution, then the problems \eqref{eq5}, \eqref{eq6} also have the unique solutions for each large $k$.
\end{itemize}

Note that the conclusion of Theorem \ref{th2} and its consequences cease to be valid for arbitrary bounded linear operators between infinite-dimensional Banach spaces.

\section{Examples}

\textit{Example 1.} Consider a linear one-point boundary-value problem for differential equation of the first order
\begin{equation}\label{1.6.1t1}
    (Ly)(t):= y'(t)+Ay(t)=f(t),\quad
t \in[a,b],
\end{equation}
\begin{equation}\label{1.3t1}
By= \sum _{k=0}^{n-1} \alpha_{k} y^{(k)}(a)=c,
\end{equation}
\noindent where $A$ is a constant $(m \times m)$ -- matrix, a vector-valued function $f(\cdot)$ belongs to the space~$(W_{p}^{n-1})^{m}$, matrices $\alpha_{k}$ belong to the space $\mathbb{C}^{l\times m}$, a vector $c \in \mathbb{C}^{l}$, $B\colon (W_{p}^{n})^{m} \rightarrow\mathbb{C}^{l}$, $(L,B)\colon (W^{n}_p)^m\rightarrow (W^{n-1}_p)^m\times\mathbb{C}^l$, $y(\cdot)\in (W_{p}^{n})^m$.

Denote by $Y(\cdot)\in (W_p^n)^{m\times m}$ the unique solution of a linear homogeneous matrix equation of the form \eqref{1.6.1t1} with the initial Cauchy condition
\begin{equation*}\label{r31}
  Y'(t)+A Y(t)=O_{m},\quad t\in (a,b), \quad Y(a)=I_{m},
  \end{equation*}
\noindent where $I_m$~is identity $(m \times m)$ -- matrix.

Put
\begin{equation*}\label{3.BY1}
M(L,B)=[BY]:=\left( B \begin{pmatrix}
                                              y_{1,1}(\cdot) \\
                                              \vdots \\
                                              y_{m,1}(\cdot) \\
                                            \end{pmatrix}
,\ldots,
                                    B \begin{pmatrix}
                                              y_{1,m}(\cdot) \\
                                              \vdots \\
                                              y_{m,m}(\cdot) \\
                                            \end{pmatrix}\right) \in \mathbb{C}^{m\times l}.
\end{equation*}
Then the fundamental matrix and its $k$-th derivative will have the following form:
\begin{gather*}
Y(t)= \operatorname{exp}\big(-A(t-a)\big), \quad Y(a) = I_{m}; \\
Y^{(k)}(t)= (-A)^k \operatorname{exp}\big(-A(t-a)\big), \quad Y^{(k)}(a) = (-A)^k, \quad k \in \mathbb{N}.
\end{gather*}

Substituting these value into the equality \eqref{1.3t1}, we have
$$M(L,B)=\sum_{k=0}^{n-1}\alpha_{k}(-A)^k.$$

Theorem \ref{th1} implies that $\operatorname{ind} \, (L, B)=\operatorname{ind} \, (M(L, B))= m-l$.

Therefore, by Theorem \ref{th2}, we obtain
\begin{gather*}
\operatorname{dim} \operatorname{ker}(L,B)=\operatorname{dim} \operatorname{ker}\left(\sum_{k=0}^{n-1}\alpha_{k}(-A)^k\right)=
m-\operatorname{rank}\left(\sum_{k=0}^{n-1}\alpha_{k}(-A)^k\right),  \\
\operatorname{dim} \operatorname{coker}(L,B)=-m+l+\operatorname{dim} \operatorname{coker}\left(\sum_{k=0}^{n-1}\alpha_{k}(-A)^k\right)=
l-\operatorname{rank}\left(\sum_{k=0}^{n-1}\alpha_{k}(-A)^k\right).\label{dimcoker}
\end{gather*}

From these formulas it follows that the Fredholm numbers of the problem do not depend on the choice of the length of the interval $(a, b)$.

\textit{Example 2.} Let us consider a multipoint boundary-value problem for the system of differential equations \eqref{1.6.1t1}, with $A(t) \equiv O_{m}$. The boundary conditions at the points $\{t_k\}_{k=0}^N \subset [a, b]$ contain derivatives of integer and$/$or \textit{fractional} orders (in the sense of Caputo \cite{Kilbas2006}). They have the next form
\begin{equation*}\label{3.BY1}
By=\sum _{k=0}^N \sum _{j=0}^s \alpha_{kj} y^{(\beta_{kj})}(t_{k})=c.
\end{equation*}
Here, numerical matrices $\alpha_{kj} \in \mathbb{C}^{l \times m}$. The nonnegative numbers $\beta_{kj}$ are such that
\begin{equation*}\label{3.BY1}
\beta_{k,0}=0 \quad \mbox{for all} \quad k \in \{1, 2, \ldots, N\}. 
\end{equation*}

Theorem \ref{th1} implies the index of the operator $(L, B)$ is equal to $m-l$.

Let us find its Fredholm numbers. In this case, the matrix $Y(\cdot)=I_{m}$. Therefore, the characteristic matrix has the form
\begin{equation*}\label{3.BY1}
M(L, B)=[BY]=\sum _{k=0}^N \sum _{j=0}^s \alpha_{kj} I_{m}^{(\beta_{kj})}=\sum _{k=0}^N \alpha_{k,0},
\end{equation*}
since the derivatives $I_{m}^{(\beta_{kj})}=0$ if $\beta_{kj} >0$. Hence, according to the Theorem \ref{th2},
\begin{gather*}
\operatorname{dim} \operatorname{ker}(L,B)=\operatorname{dim} \operatorname{ker}\left(\sum_{k=0}^{N}\alpha_{k,0}\right)=m -
\operatorname{rank}\left(\sum_{k=0}^{N}\alpha_{k,0}\right),  \\
\operatorname{dim} \operatorname{coker}(L,B)= - m+l + \operatorname{dim} \operatorname{coker}\left(\sum_{k=0}^{N}\alpha_{k,0} \right)=l -
\operatorname{rank}\left(\sum_{k=0}^{N}\alpha_{k,0} \right).
\end{gather*}

It follows from these formulas that the Fredholm numbers of the problem do not depend on the choice of the interval $(a, b)$, the points $\{t_k\}_{k=0}^N \subset [a, b]$, and the matrices $\alpha_{kj}$, with $j \geq 1$.

\textit{Example 3.} Consider a two-point boundary-value problem for a system of second-order differential equations generated by the expression
\begin{equation*}
    Ly(t):= y^{\prime \prime} (t)+Ay^{\prime}(t),\quad
t \in(a,b),
\end{equation*}
where $A$ is a constant matrix, with the boundary operator
\begin{equation*}
By=\sum^{n+1} _{k=0} \left(\alpha_{k} y^{(k)}(a)+ \beta_{k}y^{(k)}(b)\right).
\end{equation*}
Here, $\alpha_{k}$, $\beta_{k}$ are some rectangular numerical matrices. Then the operator
$$(L,B)\colon (W^{n+2}_p)^m\rightarrow (W^{n}_p)^m\times\mathbb{C}^l,$$
and the characteristic matrix
$M(L,B) \in \mathbb{C}^{2m\times l}$.

It is easy to verify that in this case
$$ Y_1(t)\equiv I_m, \quad Y_2(t)=\varphi (A, t),$$
where, for each fixed $t \in [a, b]$, the function $\varphi (\lambda, t):=1-\exp(-\lambda(t-a))\lambda^{-1}$ is an entire analytic function of the variable $\lambda \in \mathbb{C}$.

Then

$$ [BY_1]=\sum_{k=0}^{n+1} \left(\alpha_kI_m^{(k)}(a)+\beta_kI_m^{(k)}(b)\right)=(\alpha_0+\beta_0)I_m,$$
$$ [BY_2]=\sum_{k=0}^{n+1} \left(\alpha_k\varphi^{(k)}(A, a)+\beta_k\varphi^{(k)}(A, b)\right).$$

But
$$Y_2^{(k)}(t)=(-1)^kA^k\exp (-A(t-a)), \quad k \in \{0, \ldots, n+1\}.$$
Hence, we have
$$ [BY_2]=\sum_{k=0}^{n+1} \left(\alpha_kI_m+\beta_k\exp (-A(b-a))\right)(-A)^k.$$
Therefore, the characteristic block matrix
$$ M(L, B)= \left( \alpha_0+\beta_0; \sum_{k=0}^{n+1} \left(\alpha_k+\beta_k\exp (-A(b-a))\right)(-A)^k\right).$$

According to Theorem \ref{th2}, the dimensions of the kernel and cokernel of the inhomogeneous boundary-value problem are equal, respectively, to the dimensions of the kernel and cokernel of the matrix $M(L, B)$.

In particular, if $\beta_k \equiv 0$ and the problem is one-point, then the block characteristic matrix has the form
$$ M(L, B)= \left( \alpha_0; \sum_{k=0}^{n+1} \alpha_k(-A)^k\right).$$

Therefore, in this case, the Fredholm numbers of the boundary-value problem do not depend on the length of the interval $(a,b)$.

Note that the matrix $\exp (-A(b-a))$ can be found in an explicit form since every entire analytic function of a numerical matrix $A \in \mathbb{C}^{m \times m}$ is a polynomial of $A$. This polynomial is expressed via the matrix $A$ by the Lagrange--Sylvester Interpolation Formula (see, for example, \cite{Gan1959}). Its degree is no greater than $m-1$.

\section{Continuity of solutions in a parameter}\label{se4}

Let us consider parameterized by number $\varepsilon \in [0,\varepsilon_0)$, $\varepsilon_0>0$, linear boundary-value problem
\begin{equation}\label{eq7}
L(\varepsilon)y(t,\varepsilon):=y^{(r)}(t,\varepsilon) + \sum\limits_{j=1}^rA_{r-j}(t,\varepsilon)y^{(r-j)}(t,\varepsilon)=f(t,\varepsilon), \quad t\in (a,b),
\end{equation}
\begin{equation}\label{eq8}
 B(\varepsilon)y(\cdot;\varepsilon) = c(\varepsilon),
\end{equation}
where, for every fixed $\varepsilon$, matrix-valued functions $A_{r-j}(\cdot;\varepsilon) \in (W^{n}_p)^{m\times m}$, a vector-valued function $f(\cdot;\varepsilon) \in (W^{n}_p)^m$, a vector $c(\varepsilon)\in\mathbb{C}^{rm}$, $B(\varepsilon)$~is a linear continuous operator $B(\varepsilon) \colon (W^{n+r}_p)^m\rightarrow\mathbb{C}^{rm}$, and an unknown vector-valued function $y(\cdot;\varepsilon) \in (W^{n+r}_p)^m$.

It follows from Theorem \ref{th1} that the boundary-value problem \eqref{eq7}, \eqref{eq8} is a Fredholm one with index zero.

\begin{definition}\label{d2}
A solution to the boundary-value problem \eqref{eq7}, \eqref{eq8} depends continuously on the parameter $\varepsilon$ at $\varepsilon=0$ if the following two conditions are satisfied:
\begin{itemize}
\item [$(\ast)$] there exists a positive number $\varepsilon_{1}<\varepsilon_{0}$ such that, for any $\varepsilon\in[0,\varepsilon_{1})$ and arbitrary chosen right-hand sides $f(\cdot;\varepsilon)\in (W^{n}_p)^{m}$ and $c(\varepsilon)\in\mathbb{C}^{rm}$, this problem has a unique solution $y(\cdot;\varepsilon)$ that belongs to the space $(W^{n+r}_p)^{m}$;

\item [$(\ast\ast)$] the convergence of the right-hand sides $f(\cdot;\varepsilon)\to f(\cdot;0)$ in $(W_p^{n})^{m}$ and $c(\varepsilon)\to c(0)$ in $\mathbb{C}^{rm}$ implies the convergence of the solutions $y(\cdot;\varepsilon)\to y(\cdot;0)$ in $(W^{n+r}_p)^{m}$.
\end{itemize}
\end{definition}
\noindent Here and further, the limits are considered as $\varepsilon\to0+$.

Definition \ref{d2} is equivalent to the following two conditions:
\begin{itemize}
\item Operators $\big(L(\varepsilon), B(\varepsilon)\big)$ are invertible for sufficiently small $\varepsilon$;
\item $\big(L(\varepsilon), B(\varepsilon)\big)^{-1} \stackrel{s}{\longrightarrow} \big(L(0), B(0)\big)^{-1}$.
\end{itemize}

Consider the following assumptions:
\begin{itemize}
  \item [(0)] homogeneous boundary-value problem has only the trivial solution
$$L(0)y(t,0)=0,\quad t\in(a, b),\quad B(0)y(\cdot,0)=0;$$
  \item [(I)] $A_{r-j}(\cdot;\varepsilon)\to A_{r-j}(\cdot;0)$ in the space $(W^{n}_p)^{m\times m}$ for each number $j\in\{1,\ldots, r\}$;
  \item [(II)] $B(\varepsilon)y\to B(0)y$ in the space $\mathbb{C}^{rm}$ for every $y\in(W^{n+r}_p)^m$.
\end{itemize}

\begin{theorem}\label{th4} A solution to the boundary-value problem \eqref{eq7}, \eqref{eq8} depends continuously on the parameter $\varepsilon$ at $\varepsilon=0$ if and only if this problem satisfies conditions (0), (I), and (II).
\end{theorem}

This Theorem implies that if the operator $\big(L(0), B(0)\big)$ is invertible, then
$$\big(L(\varepsilon), B(\varepsilon)\big) \stackrel{s}{\longrightarrow} \big(L(0), B(0)\big) \Leftrightarrow \big(L(\varepsilon), B(\varepsilon)\big)^{-1} \stackrel{s}{\longrightarrow} \big(L(0), B(0)\big)^{-1}.$$
Note that the conclusion of Theorem \ref{th4} and its consequences cease to be valid for arbitrary bounded linear operators between infinite-dimensional Banach spaces. Note that the set of all irreversible operators is everywhere dense in the strong operator topology.

We supplement our result with a two-sided estimate of the error
$\bigl\|y(\cdot;0)-y(\cdot;\varepsilon)\bigr\|_{n+r,p}$ of the solution $y(\cdot;\varepsilon)$ via its discrepancy
$$
\widetilde{d}_{n,p}(\varepsilon):=
\bigl\|L(\varepsilon)y(\cdot;0)-f(\cdot;\varepsilon)\bigr\|_{n,p}+
\bigl\|B(\varepsilon)y(\cdot;0)-c(\varepsilon)\bigr\|_{\mathbb{C}^{rm}}.
$$
Here, we interpret $y(\cdot;0)$ as an approximate solution to the problem \eqref{eq7}, \eqref{eq8}.

\begin{theorem}\label{th5} Suppose that the boundary-value problem \eqref{eq7}, \eqref{eq8} satisfies conditions (0), (I), and (II). Then there exist positive numbers $\varepsilon_{2}<\varepsilon_{1}$ and $\gamma_{1}$, $\gamma_{2}$ such that, for any  $\varepsilon\in(0,\varepsilon_{2})$, the following two-sided estimate is true:
$$
\gamma_{1}\,\widetilde{d}_{n,p}(\varepsilon)
\leq\bigl\|y(\cdot;0)-y(\cdot;\varepsilon)\bigr\|_{n+r,p}\leq
\gamma_{2}\,\widetilde{d}_{n,p}(\varepsilon),
$$
where the quantities $\varepsilon_{2}$, $\gamma_{1}$, and $\gamma_{2}$ do not depend of $y(\cdot;\varepsilon)$ and $y(\cdot;0)$.
\end{theorem}

Thus, the error and discrepancy of the solution $y(\cdot;\varepsilon)$ to the boundary-value problem \eqref{eq7}, \eqref{eq8} are of the same degree of smallness.

\section{Comments and Remarks}
The results of Section \ref{se4} are inspired by Kiguradze's theorem in \cite{Kig1988}.

Works \cite{Atl3,Atl4} give sufficient conditions for solutions of multi-point boundary-value problems to be continuous with respect to the parameter in Sobolev spaces.

The approach used in this article can also be applied to other classes of function spaces (see, for example, \cite{MMS2016}).

In the $r=1$ case, Theorem~\ref{th1}  and Corollary~\ref{co1} are proved in~\cite{Atl1}. Theorem~\ref{th2} is also new for the systems of first order differential equations.

\section{Acknowledgements}
The research is financially supported by the grant of the Czech Academy
of Sciences, RVO:67985840.

The authors are grateful to the reviewer for careful reading of the manuscript and constructive remarks that helped improve the text.

\small Institute of Mathematics of the National Academy of Sciences of Ukraine \\
Tereshchenkivska Str. 3, 01024 Kyiv, Ukraine; \\
Institute of Mathematics of the Czech Academy of Sciences, \\
Zitna Str. 25, 115 67 Prague, Czech Republic\\
ORCID: 0000-0003-0186-3185, hatlasiuk@gmail.com

\vspace{+0.4cm}

\small Institute of Mathematics of the National Academy of Sciences of Ukraine \\
Tereshchenkivska Str. 3, 01024 Kyiv, Ukraine; \\
Institute of Mathematics of the Czech Academy of Sciences, \\
Zitna Str. 25, 115 67 Prague, Czech Republic\\
ORCID: 0000-0002-1332-1562,
vladimir.mikhailets@gmail.com

\end{document}